\documentclass[amsppt,11pt]{amsart}
\usepackage{amsmath}
\usepackage{amsfonts}
\usepackage{amssymb}
\usepackage{amscd}

\newtheorem{conjecture}{Conjecture}
\newtheorem{theorem}{Theorem}

\def\adots{\mathinner{\mkern2mu\raise0pt\hbox{.}  % antidiagonal dots
\mkern2mu\raise4pt\hbox{.}\mkern1mu
\raise7pt\vbox{\kern7pt\hbox{.}}\mkern1mu}}

\begin{document}

\bibliographystyle{ieeetr}

\title[Modular representations, old and new]
{Modular representations, old and new}\footnote{To appear in
"Buildings, Geometries and Groups", Springer Proceedings in
Mathematics (PROM)(2011)}

\author{Bhama Srinivasan}
  \address{Department of Mathematics, Statistics, and Computer Science (MC 249)\\
           University of Illinois at Chicago\\
           851 South Morgan Street\\
           Chicago, IL  60607-7045}
  \email{srinivas@uic.edu}

 \maketitle

\centerline{To the memory of Harish-Chandra}

\section{Introduction}

The art of story telling is very old in India, as is evidenced by
the great epics Ramayana and Mahabharata. Even now lectures by
scholars who tell contemporary versions of these epics are very
popular, as was a TV series on the Ramayana, where it is said that a
train stopped at a station long enough for the passengers to get out
and see the latest episode.

In this paper I would like to tell a story of modular
representations of symmetric groups, Hecke algebras and related
objects, starting from Brauer's introduction of the concept and
describing some recent developments connecting this theory with Lie
theory.

We begin with stating some conjectures in the representation theory
of finite groups, some of them being long-standing, and discuss
recent progress in them. We then discuss the ordinary, or
characteristic $0$, representation theory of finite reductive
groups, including Harish-Chandra theory and Deligne-Lusztig theory,
which also play a role in the modular representation theory. The
main theorem in the classification of blocks in the $\ell$-modular
representation theory is stated.

In recent years there has been a new direction in the modular
representation theory of symmetric groups and Hecke algebras via
connections with Lie theory, leading to the concept of graded
representation theory. A detailed, definitive account of these
developments has been given by A.Kleshchev \cite{K1}. In the second
part of the paper we give an introduction to these ideas, and hope
that the reader will then be encouraged to read \cite{K1} and other
papers in the literature.

{\bf Notation:} Given an algebra $A$, $A-mod$ (resp. $A-pmod$) is
the category of finitely generated $A$-modules (resp. projective
$A$-modules). For a category $\mathcal A$, $K_0(\mathcal A)$ is the
Grothiendieck group. For an algebra $A$, $K_0(A)$ is the
Grothiendieck group of the category $A-mod$.

For a finite group $G$, ${\rm Irr}(G)$ is the set of ordinary
irreducible characters of $G$. The usual inner product on the space
of class functions on a finite group will be denoted by $\langle -,-
\rangle$.

{\bf Acknowledgment.} I thank the organizers of the interesting and
enjoyable conference ``Buildings, Finite Geometries, Groups" in
Bangalore, August 2010, for their hospitality.

\section{Finite groups}

References for this section are (\cite{CR1}, Ch.12). All groups in
this section will be finite.

Let $G$ be a finite group. The theory of group characters goes back
to the work of Frobenius. If $\rho: G \rightarrow GL(n,K)$ is a
representation of degree $n$ over a field $K$, the character $\chi$
of $\rho$ is the class function on $G$ defined by $\chi(g)={\rm
Tr}(\rho(g))$. If $K = {\bf C}$ (or a sufficiently large field of
characteristic $0$) the (ordinary) irreducible characters form an
orthonormal basis of the space of $K$-valued class functions on $G$.
The number of ordinary characters is the number of conjugacy classes
of the group, and this leads to a "character table" for the group
which yields a lot of information about the group. We can now state

{\bf Problem 1.}

(1) Classify and describe the (ordinary) irreducible characters of
$G$ over $K$ as above, e.g. give their degrees.

(2) Describe the character table of $G$.

For some well-known groups such as the symmetric group $S_n$ or the
finite general linear group $GL(n,q)$ we can answer these questions.
In a fundamental paper in 1955 J.A.Green constructed the characters
of the groups $GL(n,q)$ (see e.g. T.A.Springer, Characters of
special groups, in \cite{B}).

Frobenius computed the first character table, that of $PSL(2,p)$. He
also defined the important notion of induced characters, which we
still use today even though there are generalizations for groups of
Lie type, for example.

The theory of modular representations, i.e. the representations of
$G$ over a field $k$ of characteristic $p$ where $p$ divides the
order of $G$  was developed by Richard Brauer starting in the
1930's. In this case  Brauer defined characters of the irreducible
representations as complex-valued functions on the $p$-regular
classes of the group, and they are now called Brauer characters.
Brauer then divided the ordinary characters into subsets called
"blocks" as follows.

Let $K$ be a sufficiently large field of characteristic $0$,
$\mathcal O$  a complete discrete valuation ring with quotient field
$K$, and $k$ a residue field of $\mathcal O$ such that the
characteristic of $k$ is $p$ . Consider the algebras $KG$,
${\mathcal O}G$, $kG$. We have ${\mathcal O}G = B_1\oplus B_2 \oplus
\ldots \oplus B_n$ where the $B_i$ are "block algebras",
indecomposable ideals of ${\mathcal O}G$. We have a corresponding
decomposition of $kG$. The principal block is the one which contains
the trivial character of $G$.

An invariant of a block $B$ of $G$ is the defect group, a $p$-
subgroup $P$ of $G$, unique up to $G$-conjugacy. One definition of
$P$ is that $P$ is minimal with respect to: Every $B$-module is a
direct summand of an induced module from $P$. The "Brauer
correspondence" then gives a bijection between blocks of $G$ of
defect group $P$ and blocks of $N_G(P)$ of defect group $P$.

An ordinary representation of $G$, i.e. a representation over $K$,
is equivalent to a representation over $\mathcal O$, and can then be
reduced mod $p$ to get a modular representation of $G$ over $k$.
Then one can define the "decomposition matrix", the transition
matrix between ordinary characters and Brauer characters. This leads
to:

\begin{itemize}
\item a partition of the ordinary characters, or $KG$-modules, into blocks
\item a partition of the Brauer characters, or $kG$-modules, into blocks
\item a partition of the decomposition matrix into blocks
\end{itemize}

Example:  $A_5$, $p=2$: Ordinary characters in the principal block:

\[\begin{array}{c|c|c|c|c|c}
\mathrm{order\ of\ element}&1&2&3&3&5\\\hline \mathrm{class\
size}&1&15&20&12&12\\\hline\hline \chi_1&1&1&1&1&1\\\hline
\chi_2&5&1&-1&0&0\\\hline
\chi_3&3&-1&0&\frac{1-\sqrt{5}}{2}-1&\frac{1+\sqrt{5}}{2}-1\\\hline
\chi_4&3&-1&0&\frac{1+\sqrt{5}}{2}-1&\frac{1-\sqrt{5}}{2}-1\\\hline

\end{array}\]

Example:  $A_5$, $p=2$: Brauer characters in the principal block:

\[\begin{array}{c|c|c|c|c}
\mathrm{order\ of\ element}&1&3&3&5\\\hline \mathrm{class\
size}&1&20&12&12\\\hline\hline \psi_1&1&1&1&1\\\hline
\psi_2&2&-1&\frac{1+\sqrt{5}}{2}-1&\frac{1-\sqrt{5}}{2}-1\\\hline
\psi_3&2&-1&\frac{1-\sqrt{5}}{2}-1&\frac{1+\sqrt{5}}{2}-1

\end{array}\]

Decomposition matrix for Principal Block of $A_5$: $\begin
{pmatrix}1&1&1&1\cr
               0&0&1&1\cr
               0&1&0&1\cr
                \end {pmatrix} $

Some of the main problems in the modular theory are:

{\bf Problem 2.}
\begin{enumerate}
\item   Describe the Brauer characters of $G$ over a (sufficiently large)
field $k$ of characteristic $p$ as above, e.g. give their degrees.
\item   Describe the blocks
\item Find the decomposition matrix $D$
\item Global to local: Describe information on the block $B$ by "local
information", i.e. from blocks of subgroups of the form $N_G(P)$,
$P$ a $p$-group
\end{enumerate}

Problem 2(1) is hard, and still open even for $S_n$ and $GL(n,q)$.
We will say more about (2) and (3) later. We discuss (4) below.

The Sylow theorems are among the first facts we learn about finite
groups. One of Brauer's insights was to realize the importance of
what we now call local subgroups of a group $G$, i.e. subgroups of
the form $N_G(P)$ where $P$ is a $p$-subgroup of $G$, for $p$-
modular representation theory. His idea is that "global information"
about a $p$-block of $G$ must be obtainable from "local
information", i.e. the characters of $N_G(P)$ where $P$ is the
defect group of the block. A striking example of this is when the
defect group of the block is cyclic (see \cite{A}, p.126).

There have been several conjectures put forward in modular
representation theory relating the characters in a block of $G$ with
the characters in blocks of local subgroups of $G$. Of these, the
simplest to state is the McKay conjecture, which asserts that for
every finite group $G$ and every prime $p$, the number of
irreducible characters of $G$ having degree prime to $p$ is equal to
the number of such characters of the normalizer of a Sylow
$p$-subgroup of $G$. There also exists a block-wise version, the
Alperin-Mckay conjecture, comparing characters in blocks of G and in
the normalizer of their defect groups. The Alperin weight conjecture
counts the number of Brauer characters (\cite{CE}, p.96;
\cite{Mar}). There are various further deep conjectures in the
representation theory of finite groups, the most elaborate ones
being the Isaacs-Navarro conjectures, Dade's conjectures and
Brou\'e's conjectures \cite{An}, \cite{Mar}.

The McKay-conjecture was proved for solvable groups and more
generally for p-solvable groups as well as for various classes of
non-solvable groups, like the symmetric groups and the general
linear groups over finite fields, but it remained open in general.
Recently, Isaacs, Malle and Navarro \cite{IMN} reduced the McKay
conjecture to a question about simple groups and gave a list of
conditions that they hoped all simple groups will satisfy. They
showed that the McKay conjecture will hold for a group $G$ if every
simple group involved in $G$ satisfied these conditions. For more on
this problem see \cite{M1} and \cite{M2}. These are conjectures at
the level of characters.

 Let $\mathcal O$ be as above. Given
a block $B$ of a group $G$ and a block $b$ of a group $H$ (e.g.
$H=N_G(D)$, $D$ defect group of $B$), a perfect isometry is a
bijection between $K_0(B)$ and $K_0(b)$ preserving certain
invariants of $B$ and $b$, and an isotypy is a collection of
compatible perfect isometries \cite{Br}. Here $B$ and $b$ are
regarded as $\mathcal O$-algebras. At the level of characters,
Brou\'e's conjecture is that there is a perfect isometry between
$K_0(B)$ and $K_0(b)$ where $H=N_G(D)$. At the module level we have
Brou\'e's abelian defect group conjecture \cite{Br}, which we state
below.

If $A$ is an ${\mathcal O}$-algebra, ${\mathcal D}^b(A)$ is the
bounded derived category of the category $A-mod$. It is a
triangulated category.

{\bf Brou\'e's Abelian Defect Group Conjecture:}

\begin{conjecture} Let $B$ be a block
of $G$ with the abelian defect group $D$, $b$ the Brauer
correspondent of $B$ in $N_G(D)$. Then ${\mathcal D}^b(B)$ and
${\mathcal D}^b(b)$ are equivalent as triangulated categories.
\end{conjecture}

A discussion of the Abelian Defect Group Conjecture is in the 1998
ICM address of J.Rickard \cite{R}. A.Marcus \cite{Mar} gives the
current status of the conjecture and describes some of the methods
used to prove it.

\section{Finite reductive groups}

References for this section are \cite{C1}, \cite{C2}, \cite{CE},
 \cite{DM}, \cite{CR2}

Until 1955 the known finite simple groups which included the
classical groups and the alternating groups were all studied
separately. A fundamental paper of Chevalley \cite{Ch} changed this
by giving a unified treatment of finite simple groups, arising from
simple Lie algebras. A treatment of this theory can be found in
\cite{C2}.

The modern view is as follows. Let $\bf G$ be a connected reductive
group defined over ${\bf F}_q$,  $F$ a Frobenius endomorphism, $F:
{\bf G} \rightarrow {\bf G}$. Then  $G= {\bf G}^F$, the group of
$F$-fixed points of $F$, is called a finite reductive group (or
finite group of Lie type). Using the
 structure of reductive algebraic
 groups we get subgroups  of ${\bf G}$ as follows. A torus is a
 closed subgroup ${\bf T} \simeq{\bf F^\times\times F^\times\times \cdots
   \times F^\times}$, where $\bf F$ is an algebraic closure of ${\bf
   F}_q$. A Levi subgroup
$\bf L$ is the  centralizer $\bf C_G(T)$  of a torus $\bf T$. A
Borel subgroup $\bf B$ is a maximal connected solvable subgroup. By
taking $F$-fixed points we get tori $T$, Borel subgroups $B$ and
Levi subgroups $L$ in $G$. In $G$ there may be several conjugacy
classes of maximal tori, but there is, up to conjugacy, a
distinguished pair $T \subset B$ of a maximal torus and a Borel
subgroup in $G$. A subgroup containing a Borel subgroup is a
parabolic subgroup of $G$. A parabolic subgroup $P$ has a Levi
decomposition $P=LV$ where $L$ is a Levi subgroup and $V$ is the
unipotent radical of $P$.

{\bf Examples.}  $G=GL(n,q), U(n,q),Sp(2n,q)$, or $SO^{\pm}(2n,q)$.
In $GL(n,q)$ Levi subgroups $L$ are isomorphic to
$\prod_iGL(m_i,q^{e_i})$.

The modular representation theory of finite reductive groups over a
field of characteristic $p$ where $p$  divides $q$ (the ``defining
charactistic") was in fact developed earlier than the ordinary
representation theory. The work of C.Curtis on this in the setting
of ``BN-pairs" is described in \cite{B}. This theory is modeled on
the highest weight theory for representations of semisimple Lie
algebras and Lie groups. Later work in this theory, again viewing
the finite groups as coming from algebraic groups via a Frobenius
morphism, is due to Andersen, Humphreys, Jantzen, Lusztig, Soergel
and others. The reader is referred to \cite{H2} for this rich theory
which we will not describe in this paper.

In the ordinary representation theory of finite reductive groups,
Curtis, Iwahori and Kilmoyer constructed the ``principal series" for
$G$ by inducing characters of a Borel subgroup and showing that the
endomorphism algebra of the induced representation is a Hecke
algebra (these days sometimes called Iwahori-Hecke algebras) (see
\cite{CR2}, Section 67). This then led to Harish-Chandra theory,
described below.

Let $\bf P$ be an $F$-stable parabolic subgroup of $\bf G$ and $\bf
L$ an $F$-stable Levi subgroup of $\bf P$ so that $L \leq P \leq G$.
Then Harish-Chandra induction is the following map: $  R^G_L:  K_0
(KL) \rightarrow K_0(KG)$, where $K$ is a sufficiently large field
of characteristic $0$ as before, such that if $\psi \in
{\text{Irr}}(L)$ then $R^G_L(\psi) =
{\text{Ind}}^{G}_{P}(\tilde\psi)$ where $\tilde\psi$ is the
character of $P$ obtained by inflating $\psi$ to $P$, using the Levi
decomposition $P=LV$ of $P$.

We say $\chi \in {\text {Irr}}(G)$ is cuspidal if $\langle \chi,
R^G_L(\psi) \rangle = 0$ for any $L \leq P< G$ where $P$ is a proper
parabolic subgroup of $G$.  The pair $(L, \theta)$ is a cuspidal
pair if $\theta \in {\text {Irr}}(L)$ is cuspidal. We then have the
main theorem of Harish-Chandra induction:

\begin{theorem} (i) Let $(L, \theta)$, $(L', \theta')$ be cuspidal
pairs. Then $\langle  R^G_L(\theta), R^G_{L'}(\theta') \rangle = 0 $
unless the pairs $(L, \theta)$, $(L', \theta')$ are $G$-conjugate.

(ii) If $\chi$ is a character of $G$, then $\langle \chi,
R^G_L(\theta) \rangle \not= 0$ for a cuspidal pair $(L, \theta)$
which is unique up to $G$-conjugacy. Thus ${\text {Irr}}(G)$ is
partitioned into Harish-Chandra families: A family is the set of
constituents of $R^G_L(\theta)$ where $(L, \theta)$ is cuspidal.
\end{theorem}

The endomorphism algebras of Harish-Chandra-induced cuspidal
representations from parabolic subgroups were then investigated and
described by Howlett and Lehrer (see \cite{C1}, Ch 10). However, not
all irreducible representations, in particular the cuspidal ones,
were obtained this way. There were some character tables as well as
the work of J.A.Green on the characters of $GL(n,q)$  which led to
the idea (attributed to I.G.Macdonald by T.A.Springer in [Cusp forms
in finite groups, \cite{B}]), that there should be families of
characters of $G$ corresponding to characters of maximal tori. This
was in fact what was proved in the spectacular paper of Deligne and
Lusztig in 1976 using $\ell$-adic cohomlogy, where $\ell$ is a prime
not dividing $q$. They introduced a map $R^{\bf G}_{\bf T}:
K_0({\overline{\bf Q}_l}T) \rightarrow K_0({\overline{\bf Q}_l} G)$
where $\bf T$ is a maximal torus of $\bf G$. Lusztig later
generalized this map, replacing maximal tori by Levi subgroups, and
we will describe this below.

Suppose $\bf L$ is an $F$-stable Levi subgroup, not necessarily in
an $F$-stable parabolic $\bf P$ of $\bf G$. Let $\ell$ be a prime
not dividing $q$. The Lusztig linear operator is a map $R^{\bf
G}_{\bf L}: K_0({\overline{\bf Q}_l}L) \rightarrow
K_0({\overline{\bf Q}_l} G)$, which has the property that every
$\chi$ in ${\text{Irr}}(G)$ is in $R^{\bf G}_{\bf T}(\theta)$ for
some $({\bf T},\theta)$, where $\bf T$ is an $F$-stable maximal
torus and   $\theta \in {\text{Irr}}(T)$.

The unipotent characters of $ G$ are the irreducible characters
$\chi$ in   $R^{\bf G}_{\bf T}(1)$ as $\bf T$ runs over $F$-stable
maximal tori of $\bf G$. Here $1$ is the trivial character of $T$.
If $L \leq P \leq G$, where $\bf P$ is a $F$-stable parabolic
subgroup, $R^{\bf G}_{\bf L}$ is just Harish-Chandra induction.

To see how the map $R^{\bf G}_{\bf L}$ is constructed using an
algebraic variety on which $G$ and $L$ act, see \cite{DM}, Ch 11.

In a series of papers and in his book \cite{L1} Lusztig classified
all the irreducible characters of $G$, provided $\bf G$ has a
connected center. (This restriction was removed by him later.) This
classification leads to two new orthonormal bases of the space $
{\mathcal C}(G)$ of ${\overline{\bf Q}_l}$-valued class functions of
$G$: the basis of ``almost characters"  and the basis of
characteristic functions of $F$-stable character sheaves on $\bf G$.
Character sheaves are certain perverse sheaves in the bounded
derived category ${\mathcal D}G$ of constructible ${\overline{\bf
Q}_l}$-sheaves on $\bf G$. Lusztig then conjectured  that the almost
characters coincide with the characteristic functions of $F$-stable
character sheaves on $\bf G$, up to a scalar multiple, if the
characteristic $p$ of $\bf F_q$ is ``almost good". This conjecture
has now been proved by T.Shoji and others in many cases, including
groups $\bf G$ with a connected center (see \cite{Sh}). If the
conjecture is true, including the precise values of certain scalars,
the character table of $G$ is determined in principle. This shows
the power of geometrical methods in representation theory. Indeed,
the theory of character sheaves is an aspect of ``geometric
representation theory", a flourishing area of research.

\section {Symmetric groups, General linear groups, Finite reductive
groups}

References for this section are \cite{J},\cite{JK}, \cite{FS1},
\cite{CE}, \cite{G1}, \cite{G2}.

Many expositions of the ordinary representation theory of the
symmetric group $S_n$ over $\bf Q$ are available (see e.g. \cite{J},
\cite{JK}). If $\chi \in {\text {Irr}}(G)$  then we can write $\chi
= {\chi}_{\lambda}$ where $\lambda$ is a partition of $n$. Then
there is a Young diagram corresponding to $\lambda$ and $p$-hooks,
$p$-cores are defined, so that we can talk of the $p$-core of a
partition $\lambda$ for a prime $p$.

In the modular representation theory, the $p$-blocks were classified
in the famous

\begin{theorem} (Brauer-Nakayama) (\cite{JK}, p.245). The characters
$\chi_\lambda$, $\chi_\mu$ of $S_n$ are in the same $p$-block ($p$
prime) if and only if $\lambda$ and $\mu$ have the same $p$-core.
\end{theorem}

However, the $p$-modular decomposition numbers are not known for
$S_n$. The modular theory has taken surprising new directions which
will be described later.

We now look at the $\ell$-modular theory for $GL(n,q)$ where $\ell$
is a prime not dividing $q$ (the ``non-defining charactistic") which
was started in \cite{FS1} and led to the development of the theory
by several authors for finite reductive groups.

Let $q$ be odd. The $\ell$-blocks were classified in \cite{FS1} for
$\ell$ odd and extended to $\ell=2$ by M.Brou\'e. Let $e$ be the
order of $q$ mod $\ell$. The unipotent characters of $GL(n,q)$ are
parametrized by partitions of $n$, and we denote the character
corresponding to a partition $\lambda$ by ${\chi}_{\lambda}$.

\medskip

\begin{theorem} (Fong-Srinivasan) ${\chi}_{\lambda}$, ${\chi}_{\mu}$
are in the same $\ell$-block of $GL(n,q)$ if and only if $\lambda$,
$\mu$ have the same $e$-core.
\end{theorem}

Example: $n=5$, $\ell$ divides $q+1$, $e=2$. Then ${\chi}_{\lambda}$
for $5,\ 32,\ 31^2,\ 2^21,\ 1^5$ are in a block. Same for $S_5$,
$p=2$.

\medskip
Example: $n=4$, $\ell$ divides $q^2+1$, $e=4$. Then
${\chi}_{\lambda}$ for $4,\ 31,\ 21^2,\  1^4$ are in a block.

$\begin {pmatrix}*&*\cr
                *&*\cr
                \end {pmatrix} $ has no $4$-hooks.

The $\ell$-blocks were then classified for various special cases,
and finally Cabanes and Enguehard (\cite{CE}, 22.9) proved the
following theorem for unipotent blocks, which are in some sense
building blocks for all $\ell$-blocks.

Definition. A unipotent block of $G$ is a block which contains
unipotent characters.

\begin{theorem} (Cabanes-Enguehard) Let $B$ be a unipotent $\ell$-block
of $G$, $\ell$ odd. Then the characters in $B$ are a union of
Lusztig families, i.e. a union of constituents of $R^G_L(\psi)$ for
various pairs $(L,\psi)$.
\end{theorem}

This theorem leads to the following

{\bf Surprise:} Brauer Theory and Lusztig Theory are compatible!

We have now described a satisfactory solution to Problem 2 (2), i.e.
the classification of $\ell$-blocks for finite reductive groups. We
now consider Problem 2 (3), the decomposition matrix. In the case of
$GL(n,q)$ the decomposition matrix is related to that of a $q$-Schur
algebra, and this will be described in Section 7. In the case of
classical groups, work has been done by Geck, Gruber, Malle and Hiss
where $\ell$ is a so-called "linear prime", using a modular version
of Harish-Chandra theory \cite{G1}, \cite{DGHM}. The problem of finding the
decomposition matrix in general is still open, even for unitary
groups. For more recent work here see \cite{G2}.

\section {Weyl groups, Cyclotomic Hecke algebras, $q$-Schur algebras}

References are \cite{Ar1}, \cite{DJ}, \cite{Ma}.

 Our story continues with the introduction of some more major
 players.

Weyl groups play an important role in the theory of finite reductive
groups. For example, the conjugacy classes of maximal tori in a
finite reductive group $G$ are parametrized by the $F$-conjugacy
classes in the Weyl group (see \cite{C1}, 3.3). There is a set of
class functions on $G$ known as ``almost characters" which involves
the characters of $W$ and plays an important role in the character
theory of $G$, as mentioned in Section 3. The ordinary characters of
Weyl groups are understood, but as mentioned above the modular
theory has open problems even for the type $A$ case, i.e. $S_n$.

The  Hecke algebras which are deformations of Weyl groups and which
arise as endomorphism algebras of induced representations of finite
reductive groups were studied by Curtis, Iwahori and Kilmoyer
(\cite{CR2}, Ch 8). These Hecke algebras and their characters are
studied from the point of view of symmetric algebras in the book by
M.Geck and G.Pffeifer \cite{GP}. With each Hecke algebra we have
certain polynomials called generic degrees, which when specialized
give the degrees of the constituents of the induced representations
mentioned above. G.Lusztig also showed that more complicated Hecke
algebras occur as endomorphism algebras of Harish-Chandra induced
cuspidal representations of Levi subgroups, and used these results
to classify the unipotent representations. Variants of these Hecke
algebras also occur in the work of Howlett and Lehrer on more
general induced representations (\cite{C1}, Ch 10).

Later Ariki-Koike and Brou\'e-Malle introduced cyclotomic Hecke
algebras which are deformations of complex reflection groups
\cite{Ar3}. These algebras also arise in the Deligne-Lusztig theory
in a mysterious way, that is, they behave as though their
specializations are endomorphism algebras of the virtual
representations $R^{\bf G}_{\bf L} $ defined above. For example, it
has been observed that there is a bijection with signs between the
constituents of $R^{\bf G}_{\bf L}(\lambda) $ for a suitable pair
$(L,\lambda)$ and characters of a complex reflection group
(\cite{BMM}, 3.2). In the paper \cite{BMM} there is developed a
theory of ``generic groups" which are defined with respect to root
data, and the Weyl group plays a central role.

The cyclotomic Hecke algebras which we are interested in are
deformations of the groups denoted $G(m,1,n)$ by Ariki. The group
$G(m,1,n)$ is isomorphic to ${\bf Z}_m \wr S_n$.

{\bf Definition} (\cite{K1}, p.420; \cite{Ar1}, 12.1).  The
cyclotomic Hecke algebra ${\mathcal H}_n$ over a field $F$ has
generators $T_1, T_2, \ldots T_{n}$ and parameters $q, v_1, v_2,
\ldots v_m$ which satisfy the relations (\cite{Ar1}, 12.1):

$(T_1-v_1)(T_1-v_2)\ldots (T_1-v_m)=0, (T_i-q)(T_i+q)=0 \ (i \geq
2),$

$T_1T_2T_1T_2=T_2T_1T_2T_1, T_iT_j=T_jT_i \ (i \geq j+2),$

$T_iT_{i-1}T_i = T_{i-1}T_iT_{i-1}\ (3 \leq i \leq n). $

In particular we get the Hecke algebra  of type $B$ or $C$ if $m=2$.
It is a classical result that the ordinary representations of
$G(m,1,n)$ are parametrized by $m$-tuples of partitions whose sizes
add up to $n$. The analogous result holds for the ordinary
representations of ${\mathcal H}_n$ in the semisimple case. Now
${\mathcal H}_n$ is semisimple if and only if  $q^iv_j-v_k \ \
(|i|<n, \ j\neq k)$ and $1+q+ \ldots q^i \ (1 \leq i <n)$ are all
non-zero (see \cite{Ar3}, 2.9). In the non-semisimple case of
interest to us where $F$ has characteristic $0$, the $v_i$ are
powers of $q$ and $q$ is an $e$-th root of unity, we have to take
$e$-restricted partitions (also called Kleshchev partitions) when
$m=1$, and Kleshchev multi-partitions when $m > 1$, to parameterize
the irreducible modules (see \cite{Ar1}, 12.1, 12.2). Here a
partition $\lambda$ is $e$-restricted if
${\lambda}_i-{\lambda}_{i+1} < e$ for all $i \geq 1$, where the
${\lambda}_i$ are the parts of $\lambda$.

In order to study the decomposition numbers for $GL(n,q)$ Dipper and
James \cite{DJ} introduced the $q$-Schur algebra. In the classical
Schur-Weyl theory the Lie group $GL(n, {\bf C})$ and $S_d$ act as
centralizers of each other on a tensor space. The Schur algebra
$S_{d,n}$ which has a basis indexed by partitions of $n$ with no
more than $d$ parts is the image of the action of $GL(n, {\bf C})$
on the tensor space. The $q$-Schur algebra  is a deformation of the
Schur algebra, defined by Dipper and James as the endomorphism
algebra of a sum of permutation modules for the Hecke algebra $H_n$
of type $A$ corresponding to $S_n$.

The $q$-Schur algebra ${\mathcal S}_q(n)$ is defined over any field
$F$ and involves an element $q \in F$. We use the same $q$ for
$H_n$.

{\bf Definition} (\cite{Ma}, p.55). ${\mathcal S}_q(n)= {\rm
End}_{H_n}{\oplus}_{\mu}M^{\mu} $.

Here $\mu$ runs over the partitions of $n$, and $M^{\mu}$ is a
permutation module for $H_n$.

At this stage we mention another important object, the affine Hecke
algebra. This has been studied for a long time, but is of special
interest here. In a paper of Ram and Ramagge \cite{RR} it is shown
how the representation theory (in characteristic 0) of Hecke
algebras of classical type can be derived from the representation
theory of the affine Hecke algebra of type $A$. Indeed, there is
also a flourishing ``Combinatorial Representation Theory", see
\cite{BR}.

\section{Lie algebras}

References are \cite{H1}, \cite{Ja} and \cite{K}.

 In the classical Cartan-Killing theory,
finite-dimensional semisimple Lie algebras are classified by Cartan
matrices and Dynkin diagrams. Let $\mathfrak{g}$ be a complex
semisimple Lie algebra. There is a root lattice and a weight lattice
associated with $\mathfrak{g}$. The root lattice has a basis of
simple roots, and there is a concept of dominant weights with
respect to this basis. The weight lattice has a basis of fundamental
dominant weights, in duality with the basis of simple roots. The
Cartan subalgebra $\mathfrak{h}$ acts on a finite-dimensional
irreducible module of $\mathfrak{g}$ by linear functions, also
called weights, of which one is a ``highest weight". The module is
then a direct sum of weight spaces

V.Kac and B.Moody generalized this theory to construct
infinite-dimensional semisimple Lie algebras which are now known as
Kac-Moody algebras. In this theory a generalized Cartan matrix $A$
is  a real matrix of the form $A=(a_{ij})$ where $a_{ii}=2$,
$a_{ij}$ are non-positive integers if $i \not= j $, and $a_{ij}= 0 $
if and only if  $a_{ji}= 0 $. Then $A$ is said to be symmetrizable
if $A=DB$ where $D$ is diagonal and $B$ is symmetric. The classical
representation theory was also extended to the Kac-Moody case, but
now we have to consider infinite-dimensional modules. In order to
get a theory resembling the finite-dimensional case one considers
modules in the category ${\mathcal O}$. Then again the module is a
direct sum of weight spaces for the Cartan subalgebra (see \cite{K},
Ch 9).

Given a generalized Cartan matrix $A$ there are (i) a semisimple Lie
algebra (Kac-Moody algebra) $\mathfrak{g}= \mathfrak{g}(A)$, (ii) a
universal enveloping algebra $U(\mathfrak{g})$, and (iii) a quantum
enveloping algebra (a Hopf algebra) $U_q(\mathfrak{g})$. The
universal enveloping algebra of a finite-dimensional semisimple Lie
algebra is a classical object, but the quantum enveloping algebra
$U_q(\mathfrak{g})$ has been studied since the 1980's and is now a
part of Lie theory. Its representation theory mirrors the classical
theory.

Given a Kac-Moody algebra $\mathfrak{g}$ we have a weight lattice
$P$ with a basis $\lbrace \Lambda_i$, $ i \in I \rbrace$ of
fundamental dominant weights and the notion of a dominant integral
weight in $P$. If $\Lambda$ is a dominant integral weight then
$\mathfrak{g}$ has an irreducible integrable module $V(\Lambda)$
which is the direct sum of weight spaces for the Cartan subalgebra
(\cite{K}, Chs 9 and 10).

The quantum enveloping algebra $U_q(\mathfrak{g})$ is generated by
$E_i, F_i, K_i^{\pm 1}$, $i \in I$ , with some relations. We have a
triangular decomposition (isomorphism of vector spaces)
$U_q(\mathfrak{g}) ={U_q}^{+}(\mathfrak{g})\otimes
{U_q}^{0}(\mathfrak{g}) \otimes {U_q}^{-}(\mathfrak{g})$ (see
\cite{Ja}, 4.3 and Theorem 4.21, \cite{Ar1}, p.21).

In the case of the quantum enveloping algebra $U_q(\mathfrak{g})$,
for $\Lambda$ as above we have again an integrable module
$V(\Lambda)$ which is the direct sum of weight spaces for
${U_q}^{0}$, the subalgebra generated by the $K_i^{\pm 1}$. Then
$V(\Lambda)$ is irreducible if $q$ is not a root of unity.

G. Lusztig proved a remarkable result regarding represenations of
quantum enveloping algebras and thus of reductive Lie algebras . He
defined a ``canonical basis" for such an algebra $\mathfrak{g}$ with
some remarkable positivity properties, which gives rise to a
distinguished basis again called a canonical basis of any
irreducible $\mathfrak{g}$-module (see \cite{Ar1}, 7.1). Kashiwara
has a similar basis known as the crystal basis.

\section{Modular Representations, New}

References are \cite{Ma}, \cite{DJ}

New modular representation theory connects decomposition numbers for
symmetric groups, cyclotomic Hecke algebras and $q$-Schur algebras
with Lie theory.

We will first discuss what is meant by ``modular representations"
for the various objects we have introduced: (1) Finite groups of Lie
type (2) Hecke algebras (3) $q$-Schur algebras.

For (1) the method of reduction mod $p$ has already been described,
and this can be applied to $S_n$ and $GL(n,q)$. In the case of $S_n$
we have the classical Specht modules which can be defined over any
field.

(2) As in (\cite{Ma}, p.133), (\cite{Ar1}, 12.1)  consider
${\mathcal H}_n $ over a field $F$. We have elements $q \neq 0$ and
$v_1, v_2, \ldots v_m$ in $F$, and a presentation for ${\mathcal
H}_n$. Then a ``Specht module" can be defined for ${\mathcal H}_n$
corresponding to a multipartition, which is irreducible over $F$ in
the semisimple case (see Section 5). In the non-semisimple case the
Specht module is no longer irreducible and we can look at its
composition factors and talk of the decomposition matrix. When $F$
has characteristic $0$ and $q$ is a root of unity Ariki's Theorem,
which will be discussed in Section 8, describes the decomposition
matrices.

(3) As in (2) a``Weyl module" can be defined for ${\mathcal
S}_q(n)$, and its composition factors when  $q$ is an  $e$-th root
of unity give rise to a decomposition matrix. Then Dipper and James
proved a remarkable relationship between the decomposition matrices
of ${\mathcal S}_q(n)$, $GL(n,q)$ and $S_n$ (\cite{Ma}, 6.47). Their
theorem states that the decomposition matrix of ${\mathcal S}_q(n)$,
now defined over a field $F$ of characteristic $\ell$, with $e$ the
order of $q$ mod $\ell$, is a submatrix of the decomposition matrix
of $FGL(n,q)$. In fact, it is the same as the part of the
decomposition matrix of $FGL(n,q)$ corresponding to the unipotent
representations. An interesting fact is that both these matrices are
square matrices; in the case of $GL(n,q)$ this fact is known by
\cite{FS1}.

However, by the work of Varagnolo-Vasserot (\cite{VV},\cite{Ma}), we
know the decomposition matrix of $\mathcal S_q(n)$ only in
characteristic $0$. One then has the James Conjecture, which is
still open, which relates the decomposition matrices in
characteristic $0$ and characteristic $\ell$ by means of an
``adjustment matrix" (\cite{Ma}, p.115).

{\bf Example, Source: GAP} An example of a decomposition matrix $D$
for ${\mathcal S}_q(n)$, $n=4$, $e=4$: $\begin
{pmatrix}4||&1&0&0&0\cr
             31||&1&1&0&0\cr
             211||&  0&1&1&0\cr
              1111||&  0&0&1&1\cr
                \end {pmatrix}$

\section{Introducing Lie Theory}

References are \cite{Ar1}, \cite{Ma}, \cite{G}.

The story now continues with a classic paper of Lascoux, Leclerc and
Thibon on the Hecke algebra $H_n$ of type $A$, defined over a field
of characteristic $0$. They gave an algorithm to compute a matrix,
conjectured to be the  decomposition matrix of $H_n$ when the
parameter $q$ is an $e$-th root of unity, and connected it with the
affine Kac-Moody algebra $\widehat {sl_e}$ (see \cite{Ma}, p.95).
Ariki proved the conjecture, in fact for cyclotomic Hecke algebras.
We state his theorem as in (\cite{Ar1}, Theorem 12.5).

Consider the cyclotomic Hecke algebra ${\mathcal H}_n$ over a field
$F$ of characteristic $0$ associated with the group $G(m,1,n)$ with
parameters $q, v_1, v_2, \ldots v_m$, where each $v_i=
q^{\gamma_i}$, $\gamma_i \in {\bf Z/eZ}, i=1,2, \ldots m$ and $q$ is
an $e$-th root of unity, not equal to $1$. Let $\Lambda = \sum_i
n_i\Lambda_i$, where $v_i$ occurs with multiplicity $n_i$. Let
$V(F)$ be the complexified Grothiendieck group $\oplus_{n
\geq0}K_0({\mathcal H}_n-{\rm pmod}) \otimes_{\bf Z}{\bf C}$.

{\bf Remark.} Here all the ${\mathcal H}_n$ are assumed to have the
same parameters. Note that $q$ is a special parameter. Then

\begin{theorem}(Ariki). The Kac-Moody Lie algebra $\widehat {sl_e}$
acts on $V(F)$ and $V(F)$ is isomorphic to the irreducible $\widehat
{sl_e}$-module of highest weight $\Lambda$. The canonical basis (in
the sense of Lusztig) of the module, specialized at $q=1$
corresponds under the isomorphism to the basis of $V(F)$  of
indecomposable projective modules of the ${\mathcal H}_n$.
\end{theorem}

Thus, we introduce Lie theory into the modular representation theory
of the Hecke algebras ${\mathcal H}_n$ taken over all $n$. Now we
also have a standard basis for $V(F)$, corresponding to the
irreducible modules of ${\mathcal H}_n$. We then have the important
fact that the transition matrix between the two bases gives the
decomposition matrices for the $e$-modular theory of the ${\mathcal
H}_n$. Furthermore, the decomposition of $V(F)$ into weight spaces
(of the Cartan subalgebra) of $\widehat {sl_e}$ corresponds to the
decomposition into blocks for the ${\mathcal H}_n$, which is another
amazing fact.

The algebra $\widehat {sl_e}$ has generators $e_i$, $f_i$, $i=1,2,
\ldots e$. Their action on $V(F)$ is given by functors called
$i$-restriction, $i$-induction which can be described
combinatorially.

It is worth noting that in the case of the symmetric groups, these
operations were defined by G.de B. Robinson, a pioneer in the
representation theory of $S_n$ (see \cite{JK}, p.271). In that case,
$i$-restriction, $i$-induction correspond to ordinary restriction
and induction followed by cutting to a block. We describe these
operations below.

\noindent Irreducible $KS_n$-modules are indexed by partitions of
$n$, and a partition of $n$ is represented by a diagram with $n$
nodes. Induction (resp. restriction) from $S_n$ to $S_{n+1}$ (resp.
$S_n \ to \ S_{n-1}$ ) is combinatorially represented by adding
(resp. removing) a node from a partition.

\noindent Now fix an integer $\ell \geq 2$. The residue $r$ at the
$(i,j)$-node of a diagram is defined as $r \equiv (j-i) (mod \
\ell)$. We define induction (resp. restriction) from $S_n$ to
$S_{n+1}$ (resp. $S_n \ to \ S_{n-1})$ by defining operators $e_i$
and $f_i$, $0 \leq i \leq (\ell -1)$, which move only nodes with
residue $i$. Then ${\rm Ind} = \sum_0^{l-1}f_i$ and ${\rm Res} =
\sum_0^{l-1}e_i$.  We have similar operations, with multipartitions,
in the case of ${\mathcal H}_n$ (see \cite{G} for a discussion).

We now see that given the  Hecke algebra ${\mathcal H}_n$ with
certain parameters, we have a dominant weight $\Lambda$ for
$\widehat {sl_e}$ associated with it. Thus the Hecke algebra can be
parametrized by $\Lambda$.  In Ariki's theorem we may denote the
Hecke algebra as ${\mathcal H}^{\Lambda}_n$ and the module $V(F)$ as
$V(\Lambda)$. This is the view adopted in \cite{K1}.

The other main idea arising from the work of Lascoux, Leclerc and
Thibon and of Ariki is the decomposition of $V(F)$ into blocks of
the ${\mathcal H}^{\Lambda}_n$ coincides with the decomposition of
the $\widehat {sl_e}$-module into weight spaces. This leads to a new
notation for the blocks:  blocks of ${\mathcal H}^{\Lambda}_n$ can
be parametrized by $\nu = \sum_{i \in I} m_ii$, where the $m_i$ are
non-negative integers and $I= {\bf Z/eZ}$ (\cite{K2}, 8.1).

\section{Categorification}

A reference is \cite{Mz}

Our story takes an abstract direction in this
section.``Categorification" is indeed one of the buzzwords of the
last few years. The idea is that set-theoretic notions are replaced
by category-theoretic notions (see \cite{Mz}). We state here a
definition given in \cite{KMS}.

 Let $A$ be a ring, $B$ a left $A$-module.
Let $a_i$ be a basis of $A$. One should find an abelian category
${\mathcal B}$ such that $K_0({\mathcal B})$ is isomorphic to $B$,
and exact endofunctors $F_i$ on $\mathcal B$ which lift the action
of the $a_i$ on $B$, i.e. the action of $[F_i]$ on $K_0({\mathcal
B})$ descends to the action of $a_i$ on $B$ so that the following
diagram commutes. The map $\phi: K_0({\mathcal B})\rightarrow B$ is
an isomorphism.

\[
\begin{CD}
K_0(\mathcal B)  @>[F_i]>>   K_0(\mathcal B) \\
@V\text{$\phi$}VV       @V\text {$\phi$}VV\\
B   @>a_i>>   B
\end{CD}
\]

{\bf Example.} Ariki's Theorem. Here $A=\widehat {sl_e}$, $ B=V(F)$
and ${\mathcal B} = \oplus_{n \geq 0}({\mathcal H}_n-{\rm pmod})$.

We return to Brou\'e's Abelian Defect Group Conjecture, stated in
Section 2. Chuang and Rouquier proved the conjecture for $S_n$ and
$GL(n,q)$ by introducing the concept of $SL_2$-categorification. On
the way they prove that if two blocks of symmetric groups (possibly
not the same $S_n$) have isomorphic defect groups, the block
algebras are derived equivalent, i.e. the derived categories of the
module categories are equivalent \cite{CR}.  This has shown the
power of categorification in a seemingly unrelated area such as
finite groups. For a discussion of categorical equivalences in the
modular representation theory of finite groups, see \cite{Mar}.

\section{KLR-algebras: the diagrammatic approach}

 References are \cite{KL}, \cite{LV}

We have seen that the algebra $\widehat {sl_e}$ is related to the
category $\sum_{n \geq0}({\mathcal H}_n-{\rm pmod})$. We have the
following question: Suppose we want to replace $U(\widehat{sl_e})$
by $U({\mathfrak g})$, where ${\mathfrak g}$ is an arbitrary
semisimple Lie algebra. Then how will we define $\mathcal B$ to give
a categorification as above?

Now we come to the papers of Khovanov-Lauda, which have connections
with Knot Theory. Given a symmetrizable Cartan matrix $A$, there is
a Kac-Moody algebra $\mathfrak{g}$ corresponding to $A$ ,  the
"negative part" of the quantum enveloping algebra
${U_q}^{-}(\mathfrak{g})$ over ${\bf Q}(q)$ and an integral form
$\mathfrak{F}$ of ${U_q}^{-}(\mathfrak{g})$, a ring defined over
${\bf Z}[q, q^{-1}]$. As before (see Section 6) $I$ is an index set
for the fundamental weights or the simple roots. Khovanov and Lauda
then construct an algebra $R$ corresponding to this data as follows
(see \cite{LV}).

A {\it diagram} is a collection of (planar) arcs connecting $m$
points on a horizontal line with $m$ points on another horizontal
line. Arcs are labeled by elements of $I$. Let ${\bf N}[I]=
\langle \nu | \nu = \sum {\nu}_i i, {\nu}_i \in {\bf N}\rangle$.
Khovanov and Lauda first define rings $R(\nu)$ for each $\nu \in
{\bf N}[I]$. The ring $R(\nu)$ is generated by diagrams in which
${\nu}_i$ arcs have the same label $i$. The product is defined by
concatenation of diagrams, subject to some relations. Then
$R=\oplus_{\nu}R(\nu)$ and $K_0(R)= \oplus_{\nu} K_0(R(\nu))$.

\begin{theorem}The
category $R-pmod
$ categorifies $\mathfrak{F}$. In other words, we
may take $A=\mathfrak{F}$, $B=\mathfrak{F}$, $\mathcal B = R-pmod$
in the commutative diagram describing the categorification. In
particular, there is an isomorphism between $K_0(R)$ and
$\mathfrak{F}$.
\end{theorem}

Rouquier defined these algebras independently and they are now
called Khovanov-Lauda-Rouquier (KLR) algebras, or sometimes as
quiver Hecke algebras. An important property of these algebras is
that they are naturally graded. This has ramifications in ``graded
representation theory" (\cite{K1}, 2.2), to be described in the next
section.

\section {Graded Representation Theory}

Ariki's theorem was regarding representations of $\widehat {sl_e}$,
but it can be generalized to those of ${U_q}(\widehat {sl_e})$
(\cite{Ar1}, 10.10). This leads to the idea of graded representation
theory, where multiplication by $q$ represents a shift in grading.
Indeed, in (\cite{K1}, p.431) the classification of graded
irreducible modules of cyclotomic Hecke algebras is stated as the
Main Problem. One can also talk of graded decomposition numbers.

We now briefly mention some of the recent important work in this
theory; references for (1) and (2) are references [38], [40] and
[42] in \cite{K1}.

(1) (Brundan and Kleshchev) Blocks of cyclotomic Hecke algebras are
isomorphic to blocks of Khovanov-Lauda-Rouquier algebras. Since the
latter are graded, there is a $\bf Z$-grading of blocks  of
cyclotomic Hecke algebras, including group algebras of symmetric
groups in positive characteristic.

(2) Building on the above, Brundan, Kleshchev and Wang have
constructed gradings of Specht modules of cyclotomic Hecke algebras.
Brundan and Kleshchev then describe graded decomposition numbers of
these Specht modules over fields of characteristic 0.

(3)  Using the connection between the Hecke algebra $H_n$ and the
$q$-Schur algebra ${\mathcal S}_q(n)$, Ariki \cite{Ar4} has shown
that ${\mathcal S}_q(n)$ and its Weyl modules can be graded.
Moreover, he shows that certain polynomials in $v$ introduced by
Leclerc and Thibon give the graded decomposition numbers for
${\mathcal S}_q(n)$ over a field of characteristic $0$ with $q$ an
$e$-th root of unity. Here we consider the algebra ${U_v}(\widehat
{sl_e})$ with a parameter $v$, as distinct from the parameter $q$ in
${\mathcal S}_q(n)$.

{\bf Example.} Graded  decomposition matrix  for $n=4$, $e=4$:
$\begin {pmatrix}4||&1&0&0&0\cr
             31||&v&1&0&0\cr
             211||&  0&v&1&0\cr
              1111||&  0&0&v&1\cr
                \end {pmatrix}$

{\bf Remark.} If $v=1$ we get the usual decomposition matrix given
in Section 7.  Other decomposition matrices are found in (\cite{Ma},
Appendix). The above matrix corresponds to a block, and the reader
can make the connection with the theorem of \cite{FS1}, stated in
Section 4; these are the partitions with an empty $4$-core.

\section{Higher Representation Theory}

A reference is \cite{Mz}.

Let $\mathcal C$ be a $2$-category. Then for any $i,j \in \mathcal
C$ the morphisms form a category whose objects are called
$1$-morphisms and morphism are called $2$-morphisms of $\mathcal C$.

 We recall the definition of Categorification in
Section 9, where the action of the algebra $A$ on $B$ was lifted to
an action by endofunctors on a category $\mathcal B$. In a good
situation the image of this lift would be a nice subcategory of the
category of endofunctors on $\mathcal B$, and we can compose these
endofunctors. Here we have an example of a $2$-category with one
object, whose $1$-morphisms are functors and $2$-morphisms are
natural transformations (see \cite{Mz}, 2.4).

The concept of ``Higher Representation Theory" was introduced by
Chuang and Rouquier \cite{CR} and continued by Rouquier and
Khovanov-Lauda. Here a 2-category is defined corresponding to a
Cartan matrix and is a categorification or a ``2-analogue" of the
enveloping algebras $U(\mathfrak{g})$ and $U_q(\mathfrak{g})$ where
$\mathfrak{g}$ is a semisimple Lie algebra. The idea then is to
study ``2-representations" of a 2-category.

One of the first explicit examples here was a categorification of
${U_q}(\widehat {sl_2})$ due to Lauda \cite{La}. B.Webster has
constructed $2$-categories corresponding to ${U_q}(\mathfrak g)$
where $\mathfrak g$ is a semisimple Lie algebra and categorifies
tensor products of irreducible representations of ${U_q}(\mathfrak
g)$. A generalization of Ariki's theorem in this context can be
found in  (\cite{W}, 5.11).

\section{End of Story?}

We have come to the end of our story. However, as we all know there
is no end to mathematical stories, and future generations will
continue them.

\end{document}